\documentclass[11pt]{article}

\usepackage{amsmath}
\usepackage{amssymb}
\usepackage{float}
\usepackage{graphicx}

\usepackage[comma,sort&compress]{natbib}
\usepackage{color} 

\usepackage{lineno}
\usepackage{url}
\usepackage{pdfpages}
\usepackage{setspace} 
\usepackage{lineno}

\usepackage{blkarray}
\usepackage{amsmath,amstext,amsthm,amssymb}
\usepackage{hyperref}

\usepackage[table]{xcolor}

\hypersetup{
	bookmarks=true,         
	unicode=false,          
	pdftoolbar=true,        
	pdfmenubar=true,        
	pdffitwindow=false,     
	pdfstartview={FitH},    
	pdftitle={My title},    
	pdfauthor={Author},     
	pdfsubject={Subject},   
	pdfcreator={Creator},   
	pdfproducer={Producer}, 
	pdfkeywords={keywords}, 
	pdfnewwindow=true,      
	colorlinks=true,       
	linkcolor=blue,          
	citecolor=blue,        
	filecolor=magenta,      
	urlcolor=cyan           
}

\usepackage[labelfont=bf,labelsep=period,justification=raggedright]{caption}

\makeatletter
\renewcommand{\@biblabel}[1]{\quad#1.}
\makeatother

\date{}

\title{\bf{Evolutionary mechanisms that promote cooperation may not promote social welfare
} }

\author{The Anh Han$^{1,\star}$, Manh Hong Duong$^{2}$, Matjaz Perc$^{3,4,5,6}$}

\begin{document}
	\maketitle
	{\footnotesize
		\noindent
		$^{1}$ School of Computing, Engineering and Digital Technologies,  Teesside University, Middlesbrough,  UK\\
  $^2$ School of Mathematics, University of Birmingham, Birmingham, UK\\
  $^3$ Faculty of Natural Sciences and Mathematics, University of Maribor, Maribor, Slovenia\\
  $^4$ Community Healthcare Center Dr. Adolf Drolc Maribor, Maribor, Slovenia\\
  $^5$ Complexity Science Hub Vienna, Vienna, Austria\\
  $^6$ Department of Physics, Kyung Hee University, Seoul, Republic of Korea\\
  
\noindent		$^\star$ Corresponding: The Anh Han (t.han@tees.ac.uk) 
	}

\newpage
\section*{Abstract}
Understanding the emergence of prosocial behaviours among
self-interested individuals is an important problem in many scientific
disciplines. Various mechanisms have been proposed to explain the
evolution of such behaviours, primarily seeking the conditions under
which a given mechanism can induce highest levels of cooperation. As
these mechanisms usually involve costs that alter individual payoffs,
it is however possible that aiming for highest levels of cooperation
might be detrimental for social welfare -- the later broadly defined
as the total population payoff, taking into account all costs involved
for inducing increased prosocial behaviours. Herein, by comparing  
stochastic evolutionary models of two well-established mechanisms of
prosocial behaviour---namely, peer and institutional incentives---we
demonstrate that the objectives of maximising
cooperation and of maximising social welfare are
often misaligned.  
 First, while peer punishment is often more effective than peer reward in promoting cooperation—especially with a higher impact-to-cost ratio—the opposite is true for social welfare. In fact, welfare typically decreases (increases) with this ratio for punishment (reward). 
Second, for institutional incentives, while maintaining similar levels of cooperation, rewards result in positive social welfare across a much broader range of parameters. Furthermore, both types of incentives often achieve optimal social welfare when their impact is moderate rather than maximal, indicating that careful planning is essential for costly institutional mechanisms to optimise social outcomes. These findings are consistent across varying mutation rates, selection intensities, and game configurations. 
 Overall, we argue for the need of adopting social welfare as
the main optimisation objective when designing and implementing
evolutionary mechanisms for social and collective goods. \\

 \noindent \textbf{Keywords:} Social welfare, cost efficiency, reward, punishment,  evolution of cooperation, social dilemma, evolutionary dynamics. 

 \newpage
 

\section{Introduction}
Since Darwin, the challenge of explaining the evolution of cooperative behavior has been actively explored across various fields, including Evolutionary Biology, Ecology, Economics, and Multi-agent Systems \citep{nowak2006,sigmund2010calculus,perc2017statistical,HanAICOM2022,paiva2018engineering,arefin2020social}. Several mechanisms have been proposed to account for the evolution of cooperation, such as kin and group selection, direct and indirect reciprocity, structured populations, pre-commitments and incentives \citep{perc2017statistical,nowak}.
Therein, the emphasis is often placed on the degree or level of cooperation that a given mechanism can induce.

However, these mechanisms typically involve costs that alter payoffs, either for the individuals involved in the interactions or for a third party (such as an institution) interested in promoting cooperation within the population. This can lead to a reduction in the overall social welfare of the population, broadly defined here as the total payoff of the population \citep{kaneko1979nash}, including all costs associated with inducing behavioral changes.
For example, let us consider peer incentives, where an agent can choose to pay a personal cost to decrease (peer punishment) or increase (peer reward) the payoff of the incentive recipient \citep{Sigmund2001PNAS,fehr2002altruistic,boyd2003evolution,han2016emergence}. Typically (and intuitively), peer punishment is considered more efficient than peer reward  as the former can lead to a higher level of cooperation, since peer punishers are more advantageous than peer rewarders when playing against defectors (see also our results in Figure \ref{fig:peerincentives_vary_impact}).  However, given that cooperative players gain an increase in payoffs when playing with rewarders, compared to no increase when playing with punishers, the overall population payoffs might be higher under peer reward even when it has a lower level of cooperation. We discuss the importance of considering social welfare for other mechanisms of prosocial behaviours and for various real-world application domains in Discussion (Section 4).

In this paper, we demonstrate that it might be more important to optimise the social welfare, rather than  focusing entirely on achieving highest  levels of cooperation. 
Because the latter can lead to a misleading, undesirable  outcome where a high cooperation level is achieved but social welfare decreases. 
We demonstrate these through analysing social welfare for two well-established classes of incentive mechanisms: peer and institutional incentives, for both positive (i.e. reward) and negative (i.e. punishment) types  \citep{Sigmund2001PNAS,van2014reward,carrotstick,duong2021cost}. 

We  adopt Evolutionary Game Theory (EGT) \citep{sigmund2010calculus,imhof2005evolutionary,smith1974theory}, a well-established mathematical framework for modelling and analysing cooperative behaviours and their emergence and stability \citep{nowak2006,perc2017statistical}. 
We derive close forms for the long-term expected social welfare, for population dynamics under varying mutation rates and selection intensities, which are key factors of Darwinian evolution \citep{nowak2006evolutionarybook}. 
Our analysis is carried out using the one-shot Prisoner's Dilemma, a well adopted game for modelling a social dilemma of cooperation \citep{coombs1973reparameterization,sigmund2010calculus,wang2015universal}.


In the next section we describe the models and methods, including derivations of social welfare and institutional costs. Results and Discussion sections will follow. We also provide additional results in the Supporting Information.

\section{Model and Methods}
\subsection{Prisoner's Dilemma}
We consider a well-mixed population where all players interact with each other via the one-shot Prisoner's Dilemma (PD) game, choosing whether to cooperate ($C$) or to defect ($D$), with payoffs given by the following payoff matrix:
\[
\begin{blockarray}{ccc}
    & C & D\\
    \begin{block}{c(cc)}
      C & R,R & S,T \\
      D & T,S & P,P \\
    \end{block}
  \end{blockarray}.
\]
If both interacting players follow the same strategy, they receive the same payoff: $R$ for mutual cooperation and $P$ for mutual defection. If the agents play different strategies, the cooperator gets the sucker's payoff $S$, and the defector gets the temptation to defect $T$. The payoff matrix corresponds to the preferences associated with the PD when the parameters satisfy the ordering $T>R>P>S$ \citep{coombs1973reparameterization}. 

 {The strength of the dilemma in the PD game can be varied adopting a simplified scaling approach from  \citep{wang2015universal,ito2018scaling,arefin2020social}. Indeed, by fixing $T - R =P - S = 1$, the dilemma strength decreases when  $R - P$ increases. }

\subsection{Evolutionary processes} 
We consider an evolutionary process of a well-mixed, finite population of $N$ interacting individuals (players). 
The players can adopt one of $m$ strategies, $1,\ ...,\ m$. 
The set of possible states of the population is 
\begin{equation}
    \Delta^m_N := \{\mathbf{n}=(n_{1},\dots,n_{m}) |\ 0 \leq n_i \leq N,  \sum_{i=1}^m n_i = N\},
\end{equation}
where  $n_i$ is the number of players currently adopting  strategy $i$ ($i=1,\ldots, m)$. In each time step of the evolutionary process, an individual $A$ is chosen at random to update their strategy. There are two ways to do so: with probability $\mu$ (mutation rate), $A$ adopts a randomly selected strategy  from the remaining $m-1$ strategies. With probability $1-\mu$, the update happens through social learning, whereby the most successful strategies tend to be imitated more often by other players (this process is equivalent to biological reproduction). That is, $A$ adopts the strategy of another, randomly chosen from the population,  player $B$, with a probability given by $P_{A,B}$. A popular approach is to use the Fermi distribution from statistical physics  {\citep{szabo1998evolutionary,traulsen2006}:  
\begin{equation}
P_{A,B} = \left(1 + e^{\beta(f_A-f_B)} \right)^{-1},
\end{equation}
}
where $f_A$ ($f_B$) denotes the fitness of individual A (B).  {The parameter $\beta$ represents the intensity of selection, i.e., how strongly the individuals base their decision to imitate on fitness comparison. For $\beta = 0$, we obtain the limit of neutral drift – the imitation decision is random. For large $\beta$, imitation becomes increasingly deterministic. Overall, 
 this} approach leads to a unified framework for evolutionary dynamics at all intensities of selection, from random drift to imitation dynamics \citep{traulsen2006,imhof2005evolutionary}.   

This elementary updating process, involving mutation and imitation, is then iterated over many time steps. As a result, we obtain an ergodic process on the space of all possible population states.
This evolutionary process defines a Markov chain with state space $\Delta^m_N$. The equilibrium of this Markov process, known as the mutation-selection distribution, is a fundamental object to quantify the evolutionary dynamics in finite populations describing the fraction of time the population spends in each population state in the long term. Understanding this equilibrium is a challenging problem due to the complexity of this calculation given the size of the transition matrix.

The number of states in the Markov chain is
\begin{equation}
S=|\Delta^m_N| = \binom{N+m  {-1}}{m  {-1}}.
\end{equation}
 {For example, in the peer incentive models we will analyse below (Section 3.1), there are three strategies, i.e. $m = 3$. Thus, in that case, the Markov chain has $S = (N+2)(N+1)/2$ states. In the institutional incentive models (Section 3.2), there are two strategies ($m = 2$) and  the Markov chain has $S = N+1$ states. }

For the transition probabilities of the Markov chain, for any two population states $\mathbf{n}$ and $\mathbf{n}^\prime$ in an evolutionary process of size $S$, the transition probability to move from $\mathbf{n}$ to $\mathbf{n}^\prime$ in one step of the process is given by 
\begin{equation}
    \label{eq: transition matrix}
\omega_{\mathbf{n}, \mathbf{n}^\prime } = \begin{cases}
\frac{n_i}{ {N}} \left(\frac{\mu}{m-1} + (1-\mu) \frac{n_j}{ {N}} P_{i,j} \right)  \quad \text{if} \ n^\prime_i = n_i -1, \ n^\prime_j = n_j +1, n^\prime_l = n_l \ \text{for } l \not\in \{i,j\}, \\
1 - \sum_{j \neq i} \frac{n_i}{ {N}} \left(\frac{\mu}{m-1} + (1-\mu) \frac{n_j}{ {N}} P_{i,j} \right) \ \  \text{if} \ \mathbf{n} = \mathbf{n}^\prime, \\
0 \ \ \text{otherwise.}
\end{cases} 
\end{equation}
 {The first equation represents the probability that a player who chooses strategy $i$ will adopt strategy $j$. It comprises the probability when it is due to mutation, i.e. $\mu \frac{n_i}{ {N}} \frac{1}{m-1}$ (the second and the third fractions represent the probability of choosing a strategist  $i$ from the population and that of selecting strategy $j$ from the set of $m-1$ strategies other than $i$), and the probability when it is due to imitation, i.e. $  (1-\mu) \frac{n_i}{ {N}}\frac{n_j}{ {N}} P_{i,j}$ (the second and the third fractions represent the probabilities of choosing a strategist  $i$ and $j$ from the population, respectively).  
}

 {Finally, b}y computing the normalised left eigenvector of the transition matrix with respect to eigenvalue $1$, we obtain the corresponding  mutation-selection (stationary) distribution.

\subsubsection{Strategy frequency}
The frequency of strategy $i$ (e.g. cooperation) is obtained by taking the average over all possible states  $\mathbf{n}$ and weighting it with the corresponding stationary distribution $\Bar{p}_{\mathbf{n}}$ 
\begin{equation}
\label{eq:frequency}
f_i = \sum_{\mathbf{n}} \frac{\mathbf{n}_i \Bar{p}_{\mathbf{n}}}{N}, \end{equation}
where $\mathbf{n}_i$ represents the quantity of individuals with strategy $i$ in state $\mathbf{n}$. 

\subsubsection{Social welfare}
Similarly, the total population payoff (social welfare), $SW$, is given as follows 
\begin{equation} 
\label{eq:SociaWelfare}
SW = \sum_{\mathbf{n}} \frac{SW(\mathbf{n}) \Bar{p}_{\mathbf{n}}}{N}, \end{equation}
where $SW(\mathbf{n})$ is the population total payoff when the population is in state $\mathbf{n}$.

\subsection{Social welfare with external  intervention}
We assume that there is an external party (i.e. an institution) that aims to promote a certain behavioural profile \citep{han2018cost,duong2021cost,wang2019exploring}. Now, when optimising social welfare one also needs to take into account the  costs spent by the third party.  

To keep it general (institutional incentives that we analyse in the present work are a special case), we assume that at state $\mathbf{n} = (n_1, \dots, n_m)$, an institutional incentive budget $\theta_\mathbf{n}$ can be used to promote a certain objective (e.g. maximising or ensure a certain threshold  of the total frequency of cooperation). 
We write $\theta_\mathbf{n} = \sum_{i=1}^m n_i\theta_i$, where $\theta_i$ is the per capita budget for strategist $i$ at state $\mathbf{n}$, which can be used for either reward or punishment.  

We denote $\Theta = \{\theta_\mathbf{n}\}_{\mathbf{n} \in \Delta^m_N}$ the overall incentive policy. 
The expected cost for providing incentive per evolutionary step is given by 
\begin{equation}
E(\Theta) = \sum_{\mathbf{n}}\theta_\mathbf{n}  \Bar{p}_{\mathbf{n}}.
\end{equation}
Thus, the total social welfare can be rewritten as follows: $$SW(\Theta) - E(\Theta).$$
Note that the population payoff depends on incentive policy $\Theta$. Namely, it alters the transition probabilities given in Equation \ref{eq: transition matrix} (more concretely, the terms $P_{i,j}$) and thus the stationary distribution.



\section{Results}
\begin{figure}
    \centering
    \includegraphics[width=0.7\linewidth]{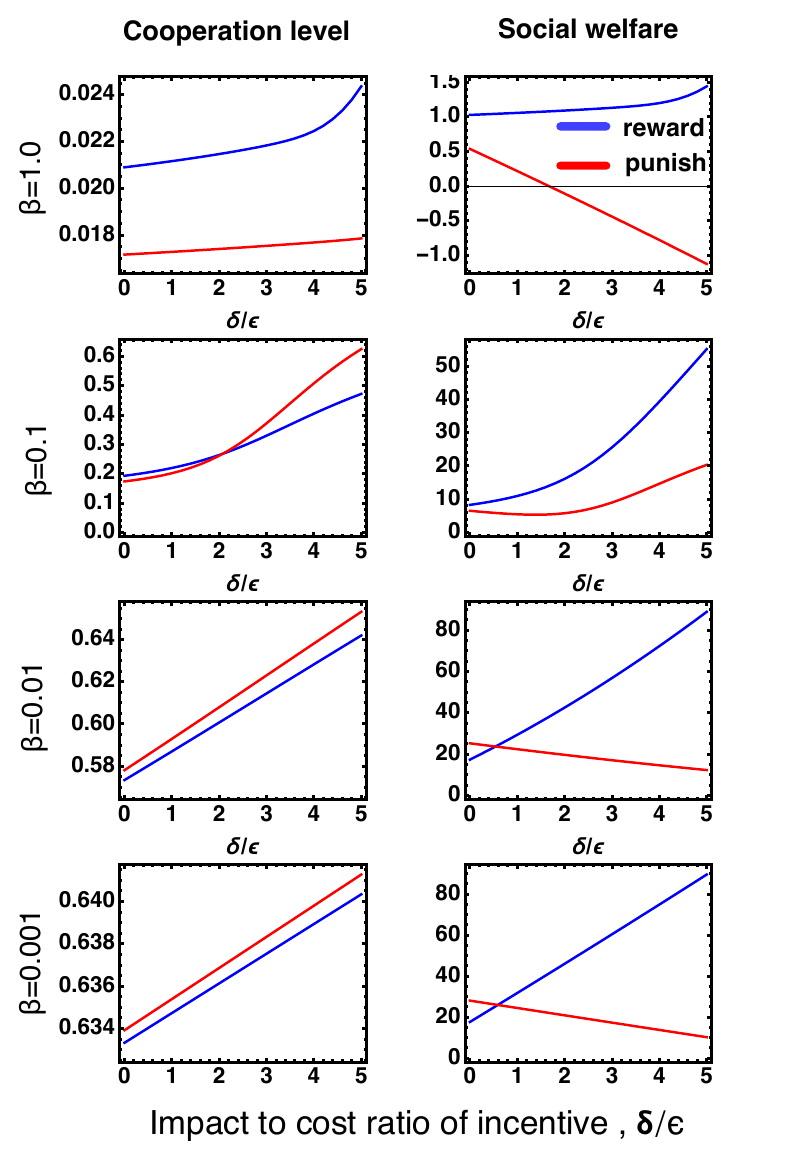}
    \caption{Impact of peer reward vs peer punishment for the long-term level of cooperation ($f_C$, see Equation \ref{eq:frequency}) and population  social welfare ($SW$, see Equation \ref{eq:SociaWelfare}), for varying the efficiency of incentive $\delta/\epsilon$, for different values of the intensities of selection $\beta$. We observe that  punishment is better than  reward for promoting cooperation in most cases, especially for weaker selection and when the impact to cost ratio of incentive is sufficiently high. However, reward leads to higher social welfare  than punishment in most cases. Parameters: Population size, $N = 50$, mutation rate $\mu = 0.01$, cost of peer incentive $\epsilon = 1$,  Prisoner's Dilemma payoff matrix $R = 1, S = -1, T = 2, P = 0$.}
    \label{fig:peerincentives_vary_impact}
\end{figure}

\subsection{Peer incentives}
A peer (social) punisher (SP) and peer (social) rewarder (SR) cooperates in the PD, and after the PD game, they  pay a cost $\epsilon$ to punish a defective  co-player or reward a cooperative one, respectively. The rewarded/punished player receives an increase/decrease of $\delta$ in their payoff.  

We consider minimal models of peer incentives in the one-shot PD game, with three strategies: unconditional cooperator (C), unconditional defector (D), and either SP or SR. The payoff matrices for peer punishment and peer reward cases  are given as follows, respectively  
\[
\begin{blockarray}{cccc}
  & C & D & SP\\
    \begin{block}{c(ccc)}
      C & R & S & R \\
      D & T & P & T -\delta \\
       SP & R & S -\epsilon & R  \\
    \end{block}
  \end{blockarray}, \ \ \ \
\begin{blockarray}{cccc}
  & C & D & SR\\
    \begin{block}{c(ccc)}
      C & R & S & R+\delta \\
      D & T & P & T  \\
       SR & R-\epsilon & S  & R -\epsilon+\delta  \\
    \end{block}
  \end{blockarray}.
\]
We now comparatively study the efficiency of the two peer incentive approaches, for promoting cooperation and enhancing the social welfare, focusing on whether these two objectives are aligned. 
Indeed, Figure \ref{fig:peerincentives_vary_impact} (left column) shows 
that, as expected,  
peer punishment usually  surpasses peer reward in effectively promoting cooperation, especially when selection is weaker and the impact to cost ratio of incentive is sufficiently high. However, interestingly, reward leads to higher social welfare  than punishment in most cases. 
 {Note that when  $\beta$ is very low, representing weak selection scenarios  (see e.g. the last row of Figure \ref{fig:peerincentives_vary_impact}, with $\beta = 0.001$), the transitions between states of the Markov chain are increasingly close  to random. That results in all strategies in the population having a similar frequency in stationary distribution. Thus, even when the impact to cost ratio is 0, the total level of cooperation (i.e. sum of frequencies of C and SP or SR) is close to $2/3$.  On the other hand, when $\beta$ is sufficiently high, representing strong selection scenarios  (see the first row of Figure \ref{fig:peerincentives_vary_impact}, with $\beta = 1.0$), the level of cooperation is very low even for the most efficient punishment/reward scenarios that are considered in our analysis (i.e. up to 5). The reason is that strong selection intensifies  the transitions that replace a C-player with a D-player among the states in the Markov chain.  Note however that if we further increase  the incentive efficiency,  high levels of cooperation can be achieved   (See Supporting Information Figure \ref{fig:peerincentives_delta_upto100}).  }

For peer reward, increasing efficiency  leads to increase in social welfare in general. It is however not  the case for punishment, where social welfare usually decreases with the impact to cost ratio. That is, applying peer punishment is often detrimental for social welfare.  

Overall, our results have shown that, in case of peer punishment,  the objective of promoting the evolution of high levels of cooperation can be detrimental for social welfare. Peer reward, on the other hand, is more efficient in promoting social welfare, even though it leads to lower levels of cooperation than punishment.  
 {The observations are robust for varying mutation rates, see Supporting Information Figure \ref{fig:peerincentives_vary_impact_different_mu}}.  {In fact, our observations indicate that a higher mutation rate further intensifies the inefficiency of peer punishment concerning social welfare.
Moreover, we  examine  the robustness of the observations for varying the  dilemma strength of the PD game, adopting a simplified scaling approach from  \citep{wang2015universal,ito2018scaling,arefin2020social}, see Supporting Information Figure \ref{fig:peerincentives_vary_impact_different_PDgames}. Again, similar observations are achieved, with  a slightly improved performance of peer punishment for a weaker dilemma strength, but only when the  intensity of selection is moderate (see  $\beta = 0.1$).}

It is noteworthy that our analysis focused on the minimal models of peer incentives, where  cooperation is generally promoted for favourable conditions of incentives (i.e. sufficiently high impact to cost ratio). These settings are suitable for the purpose of our study, as we aim to demonstrate that achieving a high level  cooperation could potentially be detrimental to social welfare. It would be interesting to examine extended models of peer incentives for example when antisocial  incentives (i.e. punishing cooperators and rewarding defectors) are included \citep{herrmann2008antisocial,rand2010anti,han2016emergence}, or when incentives are provided in a conditional or stochastic approach \citep{xiao2023evolution,chen2014probabilistic,cimpeanu2020making}.

\begin{figure}
    \centering
    \includegraphics[width=0.8\linewidth]{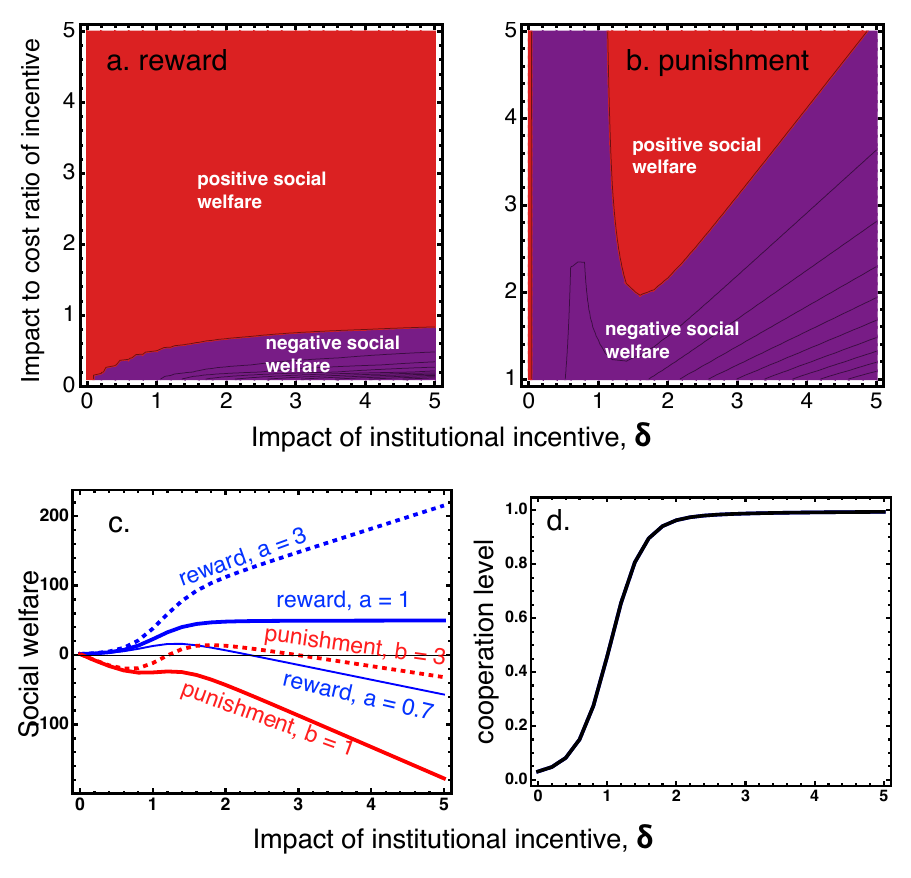}
    \caption{Impact of institutional reward vs institutional punishment for the long-term level of cooperation and population  social welfare. We observe that although both types of incentives lead to the same level of cooperation given the same incentive impact on the incentive recipient $\delta$ (assuming their impact to cost ratios are the same, i.e. $a = b$, see panel d), reward leads to positive social welfare   for a much wider range of parameters (compare red areas in panels a and b).    
    Parameters: Population size, $N = 50$, mutation rate $\mu = 0.001$, intensity of selection $\beta = 0.1$,  Prisoner's Dilemma with $R = 1, S = -1, T = 2, P = 0$.}
    \label{fig:institutional_incentive}
\end{figure}

\subsection{Institutional incentives}
We now analyse the alignment between promoting cooperation and social welfare in  minimal  models of institutional reward and punishment \citep{gois2019reward,carrotstick,duong2021cost}. Namely, we consider a well-mixed population consisting of individuals playing the one-shot PD game who can adopt either C or D in the interactions. 

The social welfare now needs to take into account the cost for providing incentives spent by the external institution  for promoting cooperation. The institution might have different levels of efficiency using the budget, which can also be different for implementing reward or punishment.  


To reward a cooperator (to punish a defector), the institution has to pay an amount $\theta/a$ ($\theta/b$, respectively) so that the cooperator's (defector's) payoff increases (decreases) by $\theta$, where $a, b > 0$ are constants representing the efficiency ratios of providing this type of incentive. 

Figure \ref{fig:institutional_incentive} shows that, as expected \citep{duong2021cost,han2022institutional,gois2019reward},  both types of incentives lead to the same level of cooperation assuming  {they} are equally effective (i.e. $a = b$) and costly. 
However, institutional reward leads to positive social welfare  (red areas in panels a and b) for a much wider range of incentive impact and cost. For reward, even when using incentive is rather cost-inefficient, e.g. when $a = 0.7$ (panel c), positive social welfare can be achieved for intermediate values of $\delta$. While for punishment, it  needs to be highly efficient ($b > 2$) for positive social welfare. Our observation are robust for different intensities of selection, see Supporting Information Figure \ref{fig:institutional_incentiv_different_betas}. 

For  cost-efficient institutional reward (i.e., $a > 1$), it is best for social welfare to maximise the impact (i.e. strong reward with a high per capital budget); for $a = 1$, the impact $\delta$ needs to reach a certain threshold and then the performance levels up. For lower $a$, there is an optimal threshold of $\delta$ for highest social welfare. 

Sufficiently strong institutional punishment (see panel c, with $b = 3$) can lead to positive social welfare, but even in this case it can be detrimental for social welfare when imposing a larger impact on  defectors. There is an optimal intermediate $\delta$ for highest social welfare. 

This is a notable observation as previous models of institutional incentives \citep{carrotstick,gois2019reward,sigmundinstitutions,sasaki2012take} do not and are unable to provide insights on how strong punishment is strong enough, focusing on promoting high levels of cooperation. In fact, previous works only consider that strong punishment is needed to ensure cooperation. 

Furthermore, we can observe that, for some value of $a$ and $b$, optimal social welfare is achieved when cooperation is not at its highest possible level (i.e. 100\% in this case). For example, for a $a = 0.7$, optimal social welfare is achieved when $\delta \approx 1.2$, and for $b = 3$, when $\delta \approx 1.6$. The corresponding levels of cooperation for those values are approximately 0.5 and 0.8. These clearly demonstrate that optimising cooperation and social welfare might not be always aligned. 

 {The above-mentioned observations are also valid when  considering varying dilemma strengths of the Prisoner's Dilemma, see Supporting Information Figure \ref{fig:institutional_incentiv_different_PDgames}. Notably, we also find that for weaker dilemmas (i.e. for larger $R - P$),  both institutional reward and punishment result in higher levels of social welfare, while as before, reward performs better than punishment for similar effect to cost ratios (i.e. $a = b$). This further reinforces the finding that the optimal institutional incentives for cooperation and social welfare are often distinct, and that increasing the costs and impacts of incentives may be counterproductive for social welfare, regardless of whether they are rewards or punishments.}












\section{Discussion}
Over the past decades, significant attention has been given to studying effective incentive mechanisms that promote the evolution of cooperation in social dilemmas \citep{nowak,sigmund2010calculus,perc2017statistical}. 
The emphasis is often placed on the extent of cooperation that a given mechanism can induce, and  the conditions regarding the parameters involved. 
Since mutual cooperation is collectively more desirable than mutation defection, ensuring high levels of cooperation usually also leads to high population welfare. 

However, as these mechanisms usually involve costs that alter individual payoffs, it is possible that aiming for highest levels of cooperation might be detrimental for social welfare. In this paper, using numerical simulations for stochastic evolutionary models for two important incentive mechanisms, peer and institutional incentives, we have demonstrated exactly that. 
 {For peer incentives, our finding showed that while peer punishment is often more effective than peer reward in promoting cooperation, in particular for a higher impact-to-cost ratio of the two incentives, the opposite is true for social welfare. In fact, the welfare typically decreases with this ratio for peer punishment, while it increases in case of peer reward. This is an important finding given that previous works usually focus more on peer punishment than peer reward, given the former advantage for enhancing cooperation \citep{Sigmund2001PNAS,herrmann2008antisocial,sigmundinstitutions,han2016emergence,cimpeanu2024digital,hauert2007}. }
 {For institutional incentives, we showed that reward fosters positive social welfare over a much wider range of parameters while maintaining similar levels of cooperation. Furthermore, both types of institutional incentives often achieve optimal social welfare when their impact is moderate rather than at the maximal level. This observation is not possible in previous models of institutional incentives as they consider only the objective of maximising cooperation \citep{gois2019reward,sasaki2012take,sigmund2010calculus}. This results  indicate that careful planning is essential for costly institutional mechanisms to optimise social goods. 
}

Closely related to the present work are a recent body of literature on cost optimisation of institutional incentives for promoting cooperation and fairness  {\citep{han2018cost,duong2021cost,DuongDurbacHan2022,wang2022incentive,cimpeanu2023social,wang2019exploring,cimpeanu2021cost,sun2021combination,cimpeanu2024digital,sun2023state}}. These works usually consider a bi-objective optimisation problem where one aims to ensure a certain minimal level of cooperation at a smallest  cost to the institution. While these works can provide insights on the optimal incentive policy, they still do not guarantee optimal social welfare. In fact, the two objectives are often misaligned (see Supporting Information, Figure \ref{fig:institutional_incentive_cost_min}).  Moreover, we argue that focusing on the social welfare provides a more convenient, single objective optimisation problem. 

We analysed here two incentive mechanisms. It might be interesting to study whether maximising social welfare  is aligned with maximising cooperation for other mechanisms of cooperation. For example, with kin selection, favouring related players might lead to reduction of social welfare due to unfair use of power to favour or patronage one's relatives or friends (aka nepotism) \citep{wilson1991nepotism}. 
For indirect reciprocity to work \citep{okada2020review}, one might need to enhance transparency of reputation, e.g. via implementing institution broadcasting reputation scores \citep{radzvilavicius2021adherence} or other costly communication mechanisms \citep{krellner2021pleasing,santos2018social}. Thus one might need to balance between promoting cooperation and the cost of enabling it.    
For direct reciprocity to work \citep{xia2023reputation}, one might need to reduce noise,  increase cognitive capacity \citep{lenaerts2024evolution,han2012corpus}, and  improve trust \citep{han2021or}, which are costly and thus need to be taken into account for balancing the  population social welfare. 
For pre-commitment to work \citep{nesse2001evolution}, it usually requires a  third party such as an institution that provides incentives \citep{han2022institutional} or broadcast commitment-based reputation  \citep{krellner2023words} for ensuring commitment compliance, which are costly and needed to be taken into account for enhanced social welfare.    
 {In the context of networked reciprocity \citep{nowak2006,xia2023reputation,santos2005scale}, different network structures and their heterogeneity can produce varying outcomes for cooperation and social welfare. Therefore, it may be necessary to consider various costly interventions, such as state-based \citep{sun2023state} or neighborhood-based incentives \citep{han2018fostering,cimpeanu2021cost}, to achieve optimal results. }

It is noteworthy that our study focused on cooperation, but the same argument would be applicable for other prosocial behaviours such as coordination, trust,  fairness, moral behaviour, technology safety development,  collective risk avoidance and pandemic intervention compliance \citep{riskSantos,capraro2021mathematical,santos2006evolutionary,cimpeanu2021cost,andras2018trusting,Han2020regulate,traulsen2023individual}.
It would be interesting to re-examine existing evolutionary mechanisms for such prosocial behaviours to see whether they promote social welfare.

Beyond prosocial behaviours such as cooperation, it is often unclear or debatable which behaviour or social norm should be promoted. In these cases, using social welfare as the optimisation objective can be particularly convenient, facilitating integrated decision-making that aims for the overall social good, especially in complex scenarios with multiple and sometimes conflicting priorities.

Such scenarios are common. For example, when designing public health programmes, determining if vaccination should be the top priority can be challenging. Utilising social welfare as the objective allows policymakers to focus on the overall health and well-being of the population.
Similarly, in education, opinions may differ on whether to prioritise STEM education, arts and humanities, vocational training, or critical thinking skills. Optimising for social welfare ensures the education system offers a balanced curriculum that supports the diverse talents and interests of students, fostering their overall development.
Another example in environment domains with debates about prioritising the reduction of carbon emissions, the preservation of biodiversity, the promotion of renewable energy, or ensuring clean water access. By focusing on social welfare, environmental policies can be created to balance these various goals, resulting in comprehensive strategies that benefit both the environment and the population.

 {We focused here on the maximisation of the  population social welfare as a whole. In many domains, it can be useful to look deeper into the distribution of the social welfare across the population, as for the same total welfare, a fairer distribution of it might be more desirable. (It is noteworthy that this is not a problem in our model as in a well-mixed population, all cooperators would receive the same  payoff.)   
For example, in the resource allocation scenarios such as those modelled by the Ultimatum and Dictator games \citep{hoffman1994preferences,zisis2015generosity,rand2013evolution,cimpeanu2021cost},  a fair distribution between resource providers and recipients are preferred, even for a lower total welfare. In addition, a critical ongoing issue in advanced  Artificial Intelligence (AI) developments is to avoid power concentration among fewer big techs \citep{verdegem2024dismantling}, in order to ensure equitable redistribution of the vast profits from advanced AI. In this domain, it would be interesting to explore mechanisms such as Windfall Clause \citep{o2020windfall}, open-source developments \citep{widder2023open}, regulatory markets \citep{clark2019regulatory,bova2024both}, and voluntary safety commitments \citep{han2022voluntary}, as potential approaches for ensuring fair and optimal social welfare.  
In all these examples, one might need to balance between efficiency (i.e. the total welfare) and equity (fair distribution of the welfare). This issue has been studied extensively in operation research literature, see e.g. a survey in \citep{karsu2015inequity}, providing suitable objective functions that combine both factors. 
}

An important ongoing initiative AI research is the design and implementation of AI for social goods \citep{tomavsev2020ai}, using AI-based solutions for effectively addressing  social problems.  To this end, there is an emerging body of evolutionary modelling studies that address prosocial behaviours in hybrid systems of human and AI agents in co-presence \citep{paiva2018engineering,santos2024prosocial,HanAICOM2022,guo2023facilitating,fernandez2022delegation,Zimmaro_Interface2024,capraro2024language, ChenInterface2024}. 
As developing AI is costly, it is crucial to understand what kind of AI are most conducive for prosocial behaviours, in a cost effective way. Thus, optimising the overall system payoff or social welfare is crucial for effective use of AI for social goods.  



\section*{Data Accessibility}
This work does not contain any data.

\section*{Declaration of AI use}
We have not used AI-assisted technologies in creating this article.

\section*{Author Contributions}
T.A.H: conceptualization, formal analysis, investigation, methodology, software, visualisation, validation, writing—original draft; 
M.H.D:  conceptualization, formal analysis, investigation, methodology, validation and writing—original draft; 
M.P: investigation, validation, and writing - review \& editing;  
All authors gave final approval for publication and agreed to be held accountable for the work performed therein.

\section*{Competing Interests}
 Authors declare no competing interests.

\section*{Acknowledgements}
T.A.H. is supported by EPSRC (grant EP/Y00857X/1) and the Future of Life Institute. M.H.D is supported by EPSRC (grant EP/Y008561/1) and a Royal International Exchange Grant IES-R3-223047. M.P. is supported by the Slovenian Research and Innovation Agency (Javna agencija za znanstvenoraziskovalno in inovacijsko dejavnost Republike Slovenije) (grants P1-0403 and N1-0232).



\newpage
 \renewcommand{\thefigure}{S\arabic{figure}}
 \renewcommand{\thetable}{S\arabic{table}}
 \setcounter{figure}{0}   
 
\section*{Supporting Information}

\begin{figure*}[ht]
    \centering
    \includegraphics[width=\linewidth]{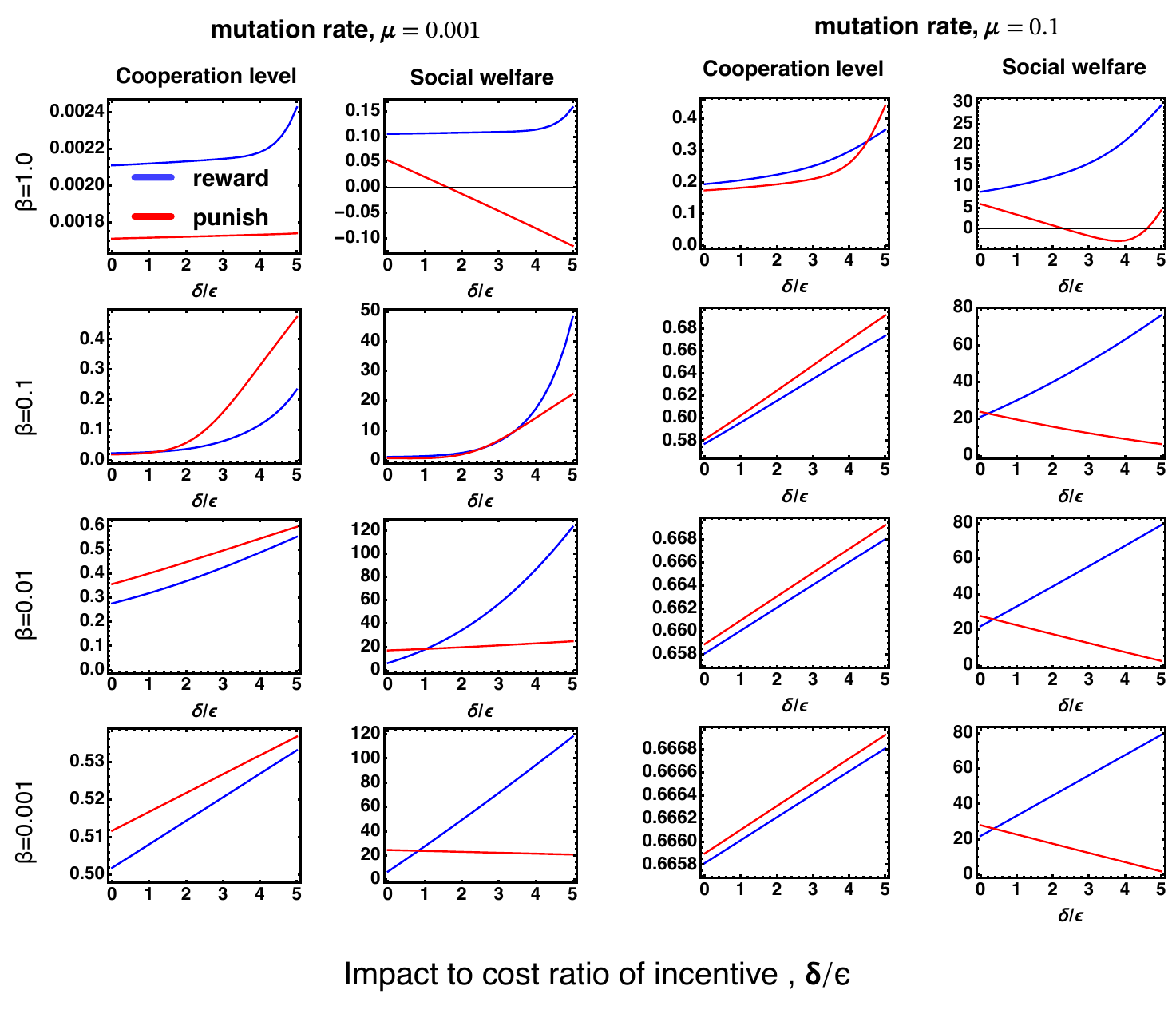}
    \caption{Impact of peer reward vs peer punishment for the long-term level of cooperation and population  social welfare, for varying the efficiency of incentive $\delta/\epsilon$, for different values of the intensities of selection $\beta$ and $\mu$. Other parameters as in Figure \ref{fig:peerincentives_vary_impact} in the main text.}
\label{fig:peerincentives_vary_impact_different_mu}
\end{figure*}

\begin{figure*}[t!]
    \centering
    \includegraphics[width=\linewidth]{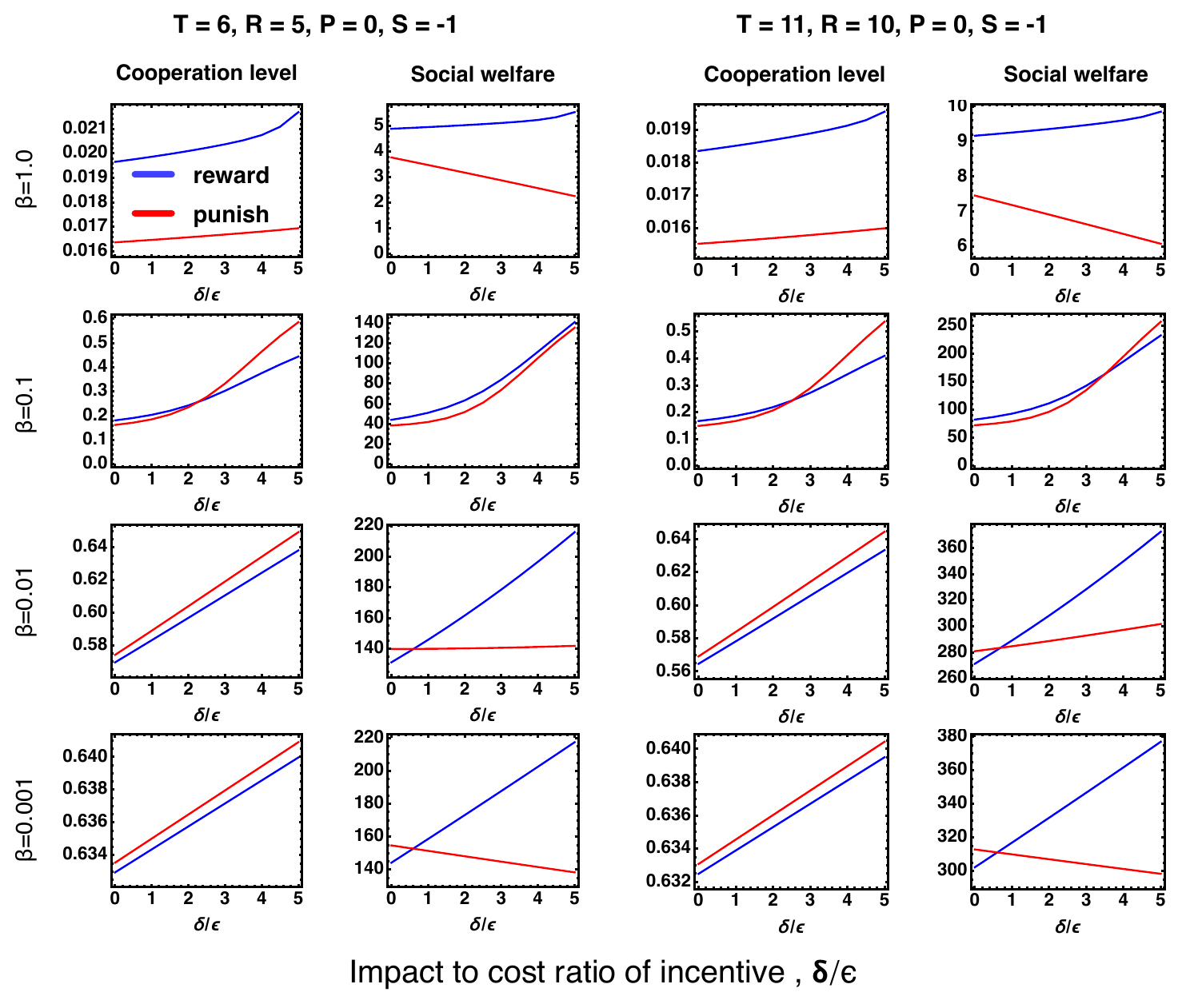}
    \caption{Impact of peer reward vs peer punishment for the long-term level of cooperation and population  social welfare, for varying the efficiency of incentive $\delta/\epsilon$, for different values of the intensities of selection $\beta$ and payoff scalings of the Prisoner's Dilemma (fixing $T - R =P - S = 1$, while varying $R - P$). Other parameters as in Figure \ref{fig:peerincentives_vary_impact} in the main text.}
\label{fig:peerincentives_vary_impact_different_PDgames}
\end{figure*}

\begin{figure*}[t!]
    \centering
    \includegraphics[width=0.8\linewidth]{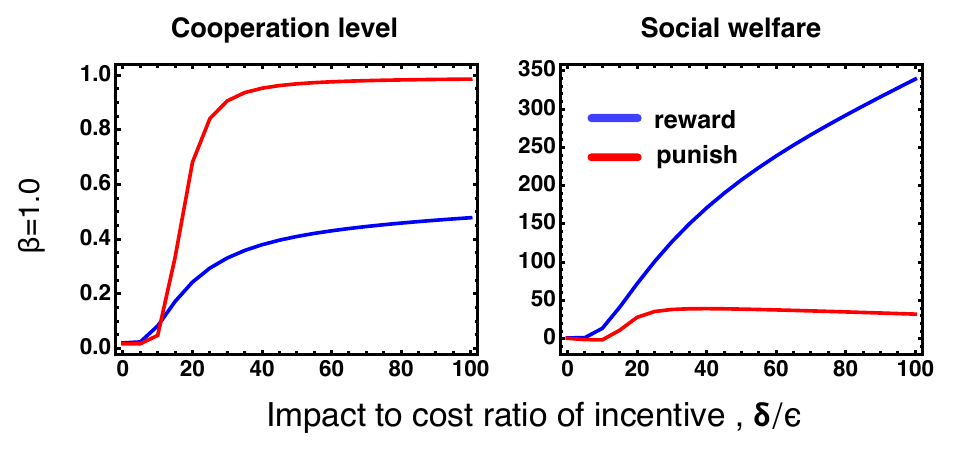}
    \caption{Impact of peer reward vs peer punishment for the long-term level of cooperation and population  social welfare, for the impact to cost ratio up to 100 under strong selection ($\beta = 1.0$). Other parameters as in Figure \ref{fig:peerincentives_vary_impact} in the main text.}
\label{fig:peerincentives_delta_upto100}
\end{figure*}

\begin{figure*}[ht]
    \centering
    \includegraphics[width=\linewidth]{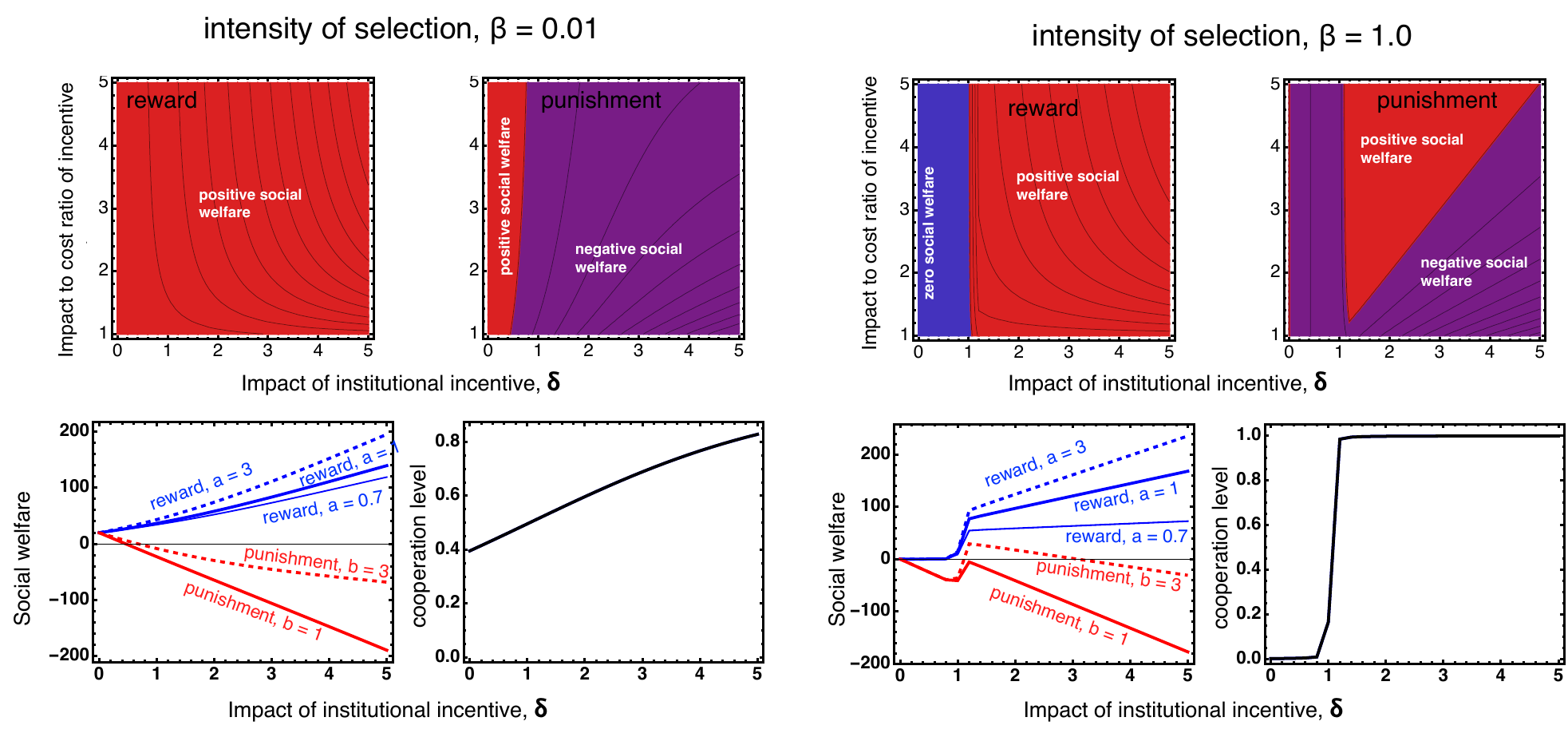}
    \caption{Impact of institutional reward vs institutional punishment for the long-term level of cooperation and population  social welfare, for different intensities of selection. Similar observations as in  Figure \ref{fig:institutional_incentive} in the main text.   
    Other parameters as in  Figure \ref{fig:institutional_incentive}.}
    \label{fig:institutional_incentiv_different_betas}
    \vspace{0.5cm}
\end{figure*}

\begin{figure*}[ht]
    \centering
    \includegraphics[width=0.8\linewidth]{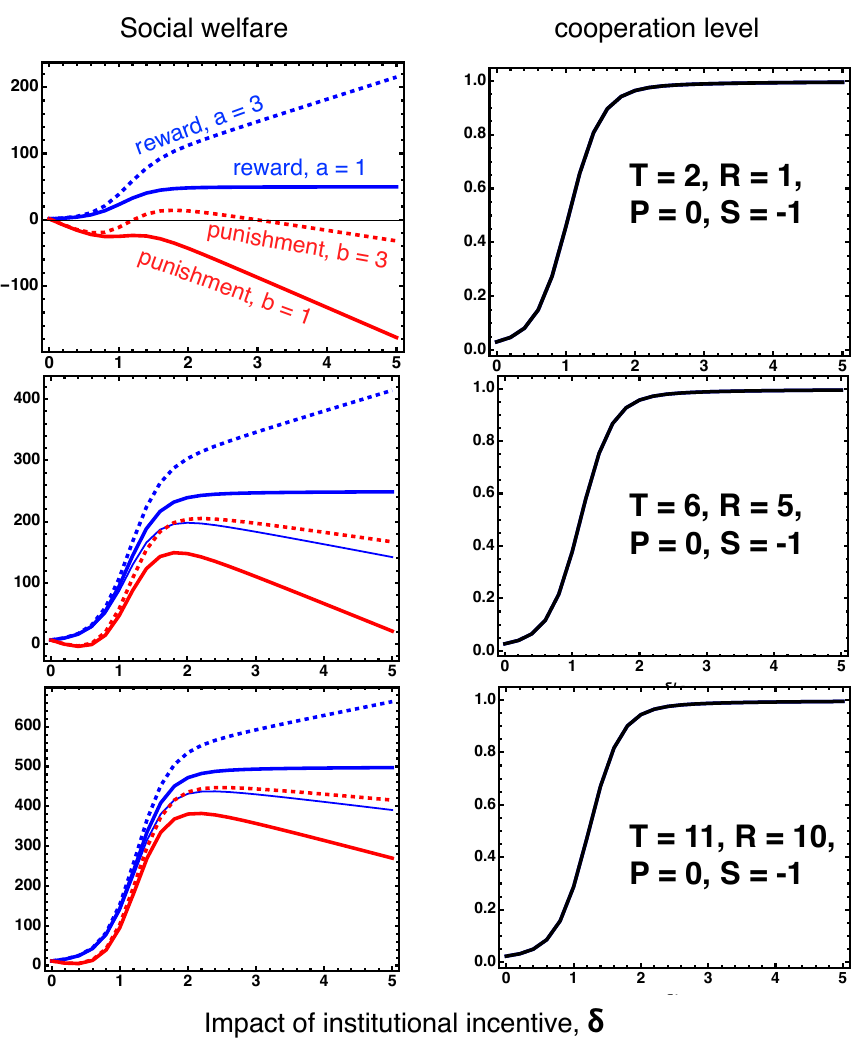}
    \caption{Impact of institutional reward vs institutional punishment for the long-term level of cooperation and population  social welfare, for different dilemma strengths of the Prisoner's Dilemma game. For weaker dilemmas, both institutional reward and punishment lead to greater levels of social welfare, and reward performs better than punishment for similar effect to cost ratios (i.e. $a = b$). 
    Other parameters as in  Figure \ref{fig:institutional_incentive}.}
    \label{fig:institutional_incentiv_different_PDgames}
    \vspace{0.5cm}
\end{figure*}

\begin{figure}
    \centering
    \includegraphics[width=\linewidth]{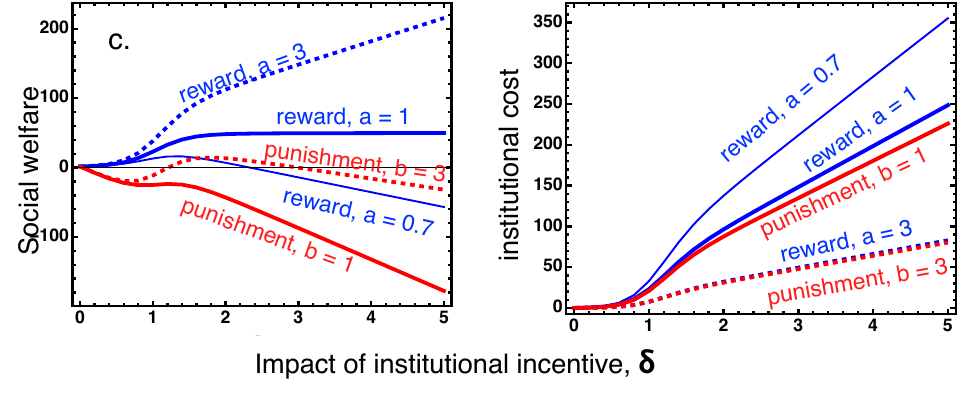}
    \caption{Maximise social welfare  vs minimise the cost of institutional incentive. 
    We observe that the these two objectives are often not aligned. For example, for highly efficient reward ($a = 3$), increasing $\delta$ leads to better social welfare. For minimising the cost while ensuring a certain desired level of cooperation (see Figure 2 in the main text), the optimal value of $\delta$ would be exactly the one that achieves the desired level of cooperation (because the institutional cost E increases with $\delta$). While for maximising social welfare, one should aim for the highest possible $\delta$ as it is an increasing function.   
    Parameters similar to Figure \ref{fig:institutional_incentive} in the main text. }
\label{fig:institutional_incentive_cost_min}
    \vspace{0.5cm}
\end{figure}

\clearpage
\newpage


\end{document}